
\magnification=1200
\pageno=1
\baselineskip=18pt
\parindent=0pt

\font\bigbf=cmbx10 scaled\magstep1
\raggedbottom

\rm
\centerline{\bigbf Local average of the hyperbolic circle problem for Fuchsian groups }
 \vskip10pt
\centerline{by Andr\'as BIR\'O\footnote{}{Research partially supported by the NKFIH (National Research, Development and InnovationOffice) Grants No. K104183, K109789, K119528, $\hbox{\rm ERC}_{
-}\hbox{\rm $\hbox{\rm HU}_{-}15\;118946$}$, and ERC-AdG Grant no. 321104}}
 \footnote{}{} \footnote{}{2000
}\hfill\break
\centerline{A. R\'enyi Institute of Mathematics, Hungarian Academy
of Sciences}

\centerline{ 1053 Budapest, Re\'altanoda u. 13-15., Hungary;
e-mail: biroand@renyi.hu}
 \vskip20pt

\noindent {\bf Abstract}. Let $\Gamma\subseteq PSL(2,{\bf R})$ be a finite volume
Fuchsian group. The hyperbolic circle problem is the
estimation of the number of elements of the $\Gamma$-orbit of
$z$ in a hyperbolic circle around $w$ of radius $R$, where $z$ and
$w$ are given points of the upper half plane and $R$ is a
large number. An estimate with error term $e^{{2\over 3}R}$ is
known, and this has not been improved for any group.
Recently Risager and Petridis proved that in the special
case $\Gamma =PSL(2,{\bf Z})$ taking $z=w$ and averaging over $z$ in a
certain way the error term can be improved to
$e^{\left({7\over {12}}+\epsilon\right)R}$. Here we show such an improvement for a
general $\Gamma$, our error term is $e^{\left({5\over 8}+\epsilon\right
)R}$ (which is better
that $e^{{2\over 3}R}$ but weaker than the estimate of
Risager and Petridis in the case $\Gamma =PSL(2,{\bf Z})$). Our main tool is
our generalization of the Selberg trace formula proved
earlier.

\vskip20pt

\noindent{\bf 1. Introduction}
\medskip

Let $H$ be the open upper half plane. The elements $\left(\matrix{
a&b\cr
c&d\cr}
\right)$
of the group $PSL(2,{\bf R})$ act on $H$ by the rule
$z\rightarrow\left(az+b\right)/\left(cz+d\right)$. Write
$$d\mu_z={{dxdy}\over {y^2}},$$
this is the $PSL(2,{\bf R})$-invariant measure on $H$.

Let $\Gamma\subseteq PSL(2,{\bf R})$ be a finite volume Fuchsian group (see
[I], p 40), i.e. $\Gamma$ acts discontinuously on $H$ and it has a fundamental domain
of finite volume (with respect to the measure $d\mu_z$). Let $F$
be a fixed fundamental domain of $\Gamma$ in $H$ (it contains exactly one point of each $
\Gamma$-equivalence class of $H$).

For $z,w\in H$ let
$$u(z,w)={{\left|z-w\right|^2}\over {4\hbox{\rm Im$z$Im$w$}}},$$
this is closely related to the hyperbolic distance $\rho (z,w)$ of $
z$
and $w$ (see [I], (1.3)). For $X>2$ define
$$N\left(z,w,X\right):=\left|\left\{\gamma\in\Gamma :\hbox{\rm \ $
4u\left(\gamma z,w\right)+2\le X$}\right\}\right|,$$
the condition here is equivalent to $\rho (z,w)\le\cosh^{-1}\left
(X/2\right)$,
hence $N\left(z,w,X\right)$ is the number of points $\gamma z$ in the
hyperbolic circle around $w$ of radius $\cosh^{-1}\left(X/2\right
).$
Therefore the estimation of $N\left(z,w,X\right)$ is called the
hyperbolic circle (or lattice point) problem. This is a
classical problem, see the Introduction of [R-P] for its
history.

In order to give the main term in the asymptotic
expansion of $N\left(z,w,X\right)$ as $X\rightarrow\infty$ we have to introduce
Maass forms.

The hyperbolic Laplace operator is given by
$$\Delta :=y^2\left({{\partial^2}\over {\partial x^2}}+{{\partial^
2}\over {\partial y^2}}\right).$$
It is well-known that $\Delta$ commutes with the action of
$PSL(2,{\bf R})$.

Let $\left\{u_j(z):\;j\ge 0\right\}$ be a complete orthonormal system of
Maass forms for $\Gamma$ (the function $u_0(z)$ is constant), let
$\Delta u_j$$=\lambda_ju_j$, where $\lambda_j=s_j(s_j-1)$, $s_j={
1\over 2}+it_j$ and $\hbox{\rm Re}s_j={1\over 2}$
or ${1\over 2}<s_j\le 1$. Note that $s_j=1$ if and only if $j=0$, and
${1\over 2}<s_j<1$ holds only for finitely many $j$.

We can now define
$$M\left(z,w,X\right):=\sqrt {\pi}\sum_{s_j\in\left({1\over 2},1\right
]}{{\Gamma\left(s_j-{1\over 2}\right)}\over {\Gamma\left(s_j+1\right
)}}u_j\left(z\right)\overline {u_j\left(w\right)}X^{s_j}.$$
It is well-known that
$$\left|N\left(z,w,X\right)-M\left(z,w,X\right)\right|=O_{z,w,\Gamma}\left
(X^{{2\over 3}}\right),$$
see e.g. [I], Theorem 12.1. The error term here has never
been improved for any group, but (as it is noted in [I]) it is conjectured that
2/3 might be lowered to any number greater than 1/2.

It was proved recently in [R-P] that in the case
$\Gamma =PSL(2,{\bf Z})$ the error term $X^{{2\over 3}}$ can be improved taking a
certain local average. More precisely, they proved that if
$f$ is a smooth nonnegative function which is compactly
supported on $F$, then
$$\int_Ff\left(z\right)\left(N\left(z,z,X\right)-M\left(z,z,X\right
)\right)d\mu_z=O_{f,\epsilon}\left(X^{{7\over {12}}+\epsilon}\right
)$$
for any $\epsilon >0$.

In the present paper we show that for this local average
the error term $X^{{2\over 3}}$ can be improved in the case of any
finite volume Fuchsian group $\Gamma$. In this generality we
get the exponent $X^{{5\over 8}+\epsilon}$, which is better than $
X^{{2\over 3}}$ but not
as strong as the result of [R-P] in the special case
$\Gamma =PSL(2,{\bf Z})$.

{\bf THEOREM 1.1.} {\it Let} $f$ {\it be a given smooth function on} $
H$ {\it such that it is compactly supported on} $F${\it , and for} $X>2$ {\it let}
$$N_f\left(X\right):=\int_Ff\left(z\right)N\left(z,z,X\right)d\mu_
z.$$
{\it Then }
$$N_f\left(X\right)=\int_Ff\left(z\right)\left(\sqrt {\pi}\sum_{s_
j\in\left({1\over 2},1\right]}{{\Gamma\left(s_j-{1\over 2}\right)}\over {
\Gamma\left(s_j+1\right)}}X^{s_j}\left|u_j\left(z\right)\right|^2\right
)d\mu_z+O_{f,\Gamma ,\epsilon}\left(X^{{5\over 8}+\epsilon}\right
)$$
{\it for every given} $\epsilon >0$.

{\bf REMARK 1.1.} As it is noted in Remark 1.3 of [R-P], the
proof there is valid only for groups similar to $PSL(2,{\bf Z})$,
as it requires strong arithmetic input not available for
most groups. Our theorem is valid for any finite volume Fuchsian group. In particular, it is valid for cocompact groups.

The main tool of our proof is our generalization of the
Selberg trace formula ([B1]), which is valid for every
finite volume Fuchsian group $\Gamma$.  Note that the operator whose trace is studied in [B1] has been used and
analysed also in a series of papers by Zelditch (see [Z1], [Z2], [Z3]).

{\bf REMARK 1.2.} As a very brief indication of the idea of our
proof we mention that
$$N\left(z,z,X\right)=\sum_{\gamma\in\Gamma}k\left(u(z,\gamma z)\right
),\eqno (1.1)$$
where $k$ is the characteristic function of the interval
$[0,x]$ with a large real $x$ (in fact $x=(X-2)/4$). We use
the decomposition
$$k\left(v\right)=k^{\ast}\left(v\right)+\left(k\left(v\right)-k^{
\ast}\left(v\right)\right),\eqno (1.2)$$
where $k^{\ast}$ is a certain smoothed version of $k$. We will estimate
the contribution of $k^{\ast}$ in (1.1) in the traditional way,
using the spectral expansion of the automorphic kernel
function given by $k^{\ast}$ and estimating the
Selberg-Harish-Chandra transform of $k^{\ast}$. However, the
contribution of $k\left(v\right)-k^{\ast}\left(v\right)$ to $N_f\left
(X\right)$ is estimated in a
completely different way, using our generalization of the
Selberg trace formula ([B1]).

{\bf REMARK 1.3.} We have seen that $N\left(z,z,X\right)$ is the number of points
in the $\Gamma$-orbit of $z$ in a hyperbolic circle around $z$ of
large radius. Note that the analogous quantity in the
euclidean case (if we choose in place of $\Gamma$ the group of translations on
the euclidean plane with vectors having integer
coordinates in place of $\Gamma$) is independent of $z$, hence
averaging in $z$ does not help in the euclidean case, the
problem there remains the same.

We mention that different kind of averages were considered by Chamizo in [C]. In particular, he proved a strong estimate for the integral with respect to $z$ of the square of $N\left(z,w,X\right)$ for a fixed $w$ over the whole fundamental domain in the case of a cocompact group, see Corollary 2.2.1 of [C].

{\bf REMARK 1.4.} The structure of the paper is the following.
In Section 2 we introduce the necessary notations. In
Section 3 we give the two main lemmas (Lemmas 3.3 and
3.4) needed for the proof of the theorem. Lemma 3.3 is
our main new tool, this is a consequence of our
generalization of the Selberg trace formula in [B1].
Lemma 3.4 is the well-known spectral expansion of an
automorphic kernel function. The proof of Theorem 1.1
is given in Section 4, using some results proved only
later on special functions and automorphic functions in Sections 5 and 6, respectively.

\noindent{\bf 2. Further notations}
\medskip

We fix a complete set $A$ of inequivalent cusps of $\Gamma$, and
we will denote its elements by $a$, $b$ or $c$, so e.g. $\sum_a$$
\sum_c$ or $\cup_a$ will
mean that $a$ and $c$ run over $A$. We say that $\sigma_a$ is a {\it scaling }
{\it matrix\/} of a cusp $a$ if $\sigma_a\infty =a$, $\sigma_a^{-
1}\Gamma_a\sigma_a=B$, where $\Gamma_a$ is
the stability group of $a$ in $\Gamma$, and $B$ is the group of
integer translations. The scaling matrix is determined up
to composition with a translation from the right.

We also fix a complete set $P$ of representatives of
$\Gamma$-equivalence classes of the set
$$\left\{z\in H:\hbox{\rm \ $\gamma z=z\hbox{\rm \ for some $id\neq$}$}
\gamma\in\Gamma\right\}.$$
$ $For a $p\in P$ let $m_p$ be the order of the stability group of $
p$ in $\Gamma$.

Let
$$P\left(Y\right)=\left\{z=x+iy:\;0<x\le 1,\>y>Y\right\},$$
and let $Y_{\Gamma}$ be a constant (depending only on the group $
\Gamma$)
such that for any $Y\ge Y_{\Gamma}$ the cuspidal zones
$F_a\left(Y\right)=\sigma_aP\left(Y\right)$ are disjoint, and the fixed fundamental
domain $F$ of $\Gamma$ is partitioned into
$$F=F\left(Y\right)\cup\bigcup_aF_a\left(Y\right),$$
where $F\left(Y\right)$ is the central part,
$$F\left(Y\right)=F\setminus\bigcup_aF_a\left(Y\right),$$
and $F\left(Y\right)$ has compact closure.

For $j\ge 0$ and $a\in A$ we have the Fourier expansion
$$u_j\left(\sigma_az\right)=\beta_{a,j}\left(0\right)y^{1-s_j}+\sum_{
n\neq 0}\beta_{a,j}\left(n\right)W_{s_j}\left(nz\right),$$
where $W$ is the Whittaker function.

The Fourier expansion of the Eisenstein series (as in [I],
(3.20)) is given by
$$E_c\left(\sigma_az,s\right)=\delta_{ca}y^s+\phi_{c,a}\left(s\right
)y^{1-s}+\sum_{n\neq 0}\phi_{a,c}\left(n,s\right)W_s\left(nz\right
).$$

For $Y\ge Y_{\Gamma}$ let us define the truncated Eisenstein series
(as in [I], pp 95-96) in the following way: for a given
$c\in A$ and every $a\in A$ let
$$E_c^Y\left(z,s\right)=E_c\left(z,s\right)-\left(\delta_{ca}\left
(\hbox{\rm Im}\sigma_a^{-1}z\right)^s+\phi_{c,a}\left(s\right)\left
(\hbox{\rm Im}\sigma_a^{-1}z\right)^{1-s}\right)\hbox{\rm \ for $ $}
z\in F_a\left(Y\right),$$
let
$$E_c^Y\left(z,s\right)=E_c\left(z,s\right)\hbox{\rm \ for $ $}z\in
F\left(Y\right),$$
finally let $E_c^Y\left(\gamma z,s\right)=E_c^Y\left(z,s\right)$ for $
\gamma\in\Gamma$ and $z\in F$.

For $Y\ge Y_{\Gamma}$ and $j\ge 0$ let us also define the truncation of
$u_j$ in the following way: for every $a\in A$ let
$$u_j^Y\left(z\right)=u_j\left(z\right)-\beta_{a,j}\left(0\right)\left
(\hbox{\rm Im}\sigma_a^{-1}z\right)^{1-s_j}\hbox{\rm \ for $ $}z\in
F_a\left(Y\right),$$
let
$$u_j^Y\left(z\right)=u_j\left(z\right)\hbox{\rm \ for $ $}z\in F\left
(Y\right),$$
finally let $u_j^Y\left(\gamma z\right)=u_j^Y\left(z\right)$ for $
\gamma\in\Gamma$ and $z\in F$.

Let $\left\{S_l:\;l\in L\right\}$ be the set of the poles in the half-plane
$\hbox{\rm Re$s>{1\over 2}$}$ of the Eisenstein series for $\Gamma$. Then ${
1\over 2}<S_l\le 1$ for every $l\in L$, and $L$ is a
finite set. We have $\beta_{a,j}\left(0\right)=0$ if $u_j\left(z\right
)$ is not
a linear combination of the residues of Eisenstein series,
so if $j\ge 0$ is such that $\beta_{a,j}\left(0\right)\neq 0$ for some $
a$, then
$s_j=S_l$ for some $l\in L$. In particular, $u_j^Y$ is the same as
$u_j$ for all but finitely many $j.$

The constants in the symbols $O$ will depend on the group
$\Gamma$. For a function $g$ we will denote
its $j$th derivative by $g^{(j)}$.

For $\lambda\le 0$ define the special function $f_{\lambda}(\theta
)$ in the
following way: $f_{\lambda}(\theta )$ is the unique even solution of the
differential equation
$$f^{(2)}(\theta )={{\lambda}\over {\cos^2\theta}}f(\theta ),\qquad
\theta\in (-{{\pi}\over 2},{{\pi}\over 2})\eqno (2.1)$$
with $f_{\lambda}(0)=1$. Note that this differential equation (which
appeared in [B1] and also in [Hu], equations (10-(11)) is the
Laplacian on functions depending only on the hyperbolic distance
from the imaginary real axis, i.e. if for $z\in H$ we write
$z=re^{i\left({{\pi}\over 2}+\theta\right)}$ with $r>0$ and $\theta
\in (-{{\pi}\over 2},{{\pi}\over 2})$, then for
$$F\left(z\right):=f_{\lambda}(\theta )$$
we have $\Delta F$$=\lambda F$.

For $\lambda\le 0$ define also the special function $g_{\lambda}(
r)$ ($r\in [0,\infty )$)
as the unique solution of
$$g^{(2)}(r)+{{\cosh r}\over {\sinh r}}g^{(1)}(r)=\lambda g(r)$$
with $g_{\lambda}(0)=1.$ Note that it is well-known (see e.g. [I], (1.20)) that this differential equation is the
Laplacian on functions depending only on the hyperbolic
distance $\rho (z,i)$ from the given point $i$, i.e. if for $z\in
H$ we write
$$G\left(z\right):=g_{\lambda}\left(\rho (z,i)\right),$$
then we have $\Delta G$$=\lambda G$.

Note that $f_0(\theta )$ and $g_0(r)$ are the identically $1$ functions.

If $m$ is a compactly supported continuous function on
$[0,\infty )$, let (see [I], (1.62))
$$g_m\left(a\right)=2q_m\left({{e^a+e^{-a}-2}\over 4}\right),\hbox{\rm \ where $
q_m$$\left(v\right)=\int_0^{\infty}{{m\left(v+\tau\right)}\over {\sqrt {
\tau}}}d\tau$},\eqno (2.2)$$
and let
$$\hbox{\rm $h_m\left(r\right)=\int_{-\infty}^{\infty}g_m\left(a\right
)e^{ira}da.$}\eqno (2.3)$$

For $\gamma\in\Gamma$ denote by $[\gamma ]$ the conjugacy class of $
\gamma$ in $\Gamma$,
i.e.
$$[\gamma ]=\left\{\tau^{-1}\gamma\tau :\,\tau\in\Gamma\right\}.$$

We use the general notation
$$\left(F,G\right)=\int_FF\left(z\right)\overline {G\left(z\ \right
)}d\mu_z.$$

We write $ $$F\left(\matrix{\alpha ,\beta\cr \gamma\cr} ;z\right)$
for the Gauss hypergeometric function. We use the notations $\Gamma\left(X\pm Y\right)=\Gamma\left(X+Y\right
)\Gamma\left(X-Y\right)$.

\noindent{\bf 3. Basic lemmas}
\medskip

Our two main results here are Lemmas 3.3 and 3.4, but
first we have to prove two simple lemmas.

{\bf LEMMA 3.1.} {\it Let} $a\in A${\it . If} $Y$ {\it is large enough (depending only on} $
\Gamma${\it ), then for} $z\in F_a\left(Y\right)$ {\it and} $\gamma
\in\Gamma$ {\it we have either}
$$u\left(\gamma z,z\right)\ge D_{\Gamma}Y^2\eqno (3.1)$$
{\it with some constant} $D_{\Gamma}>0$ {\it depending only on} $
\Gamma${\it , or we have} $\gamma\in\Gamma_a${\it .}

{\it Proof.\/} Let $z\in F_a\left(Y\right)=\sigma_aP\left(Y\right
)$, then $z=\sigma_aw$ with some
$w\in P\left(Y\right)$, and for $\gamma\in\Gamma$ we have
$$u\left(\gamma z,z\right)=u\left(\sigma_a^{-1}\gamma\sigma_aw,w\right
).\eqno (3.2)$$
Let $\sigma_a^{-1}\gamma\sigma_a=\left(\matrix{\ast&\ast\cr
C&D\cr}
\right)$. Assume $\left|C\right|>0$, then
$\hbox{\rm Im}\sigma_a^{-1}\gamma\sigma_aw={{\hbox{\rm Im$w$}}\over {\left
|Cw+D\right|^2}}\le{1\over {C^2\hbox{\rm Im$w$}}}$. Since $\hbox{\rm Im}
w>Y$, for large
enough $Y$ this implies
$$u\left(\sigma_a^{-1}\gamma\sigma_aw,w\right)\ge{{\left|Y-{1\over {
C^2Y}}\right|^2}\over {4Y{1\over {C^2Y}}}}.\eqno (3.3)$$
Since
$$\min\left\{\left|C\right|>0:\hbox{\rm \ $\left(\matrix{\ast&\ast\cr
C&\ast\cr}
\right)\in\sigma_a^{-1}\Gamma\sigma_a$}\right\}$$
exists (see [I], p 53), (3.2) and (3.3) imply (3.1).

Assume $C=0$. Then $\sigma_a^{-1}\gamma\sigma_a\infty =\infty$, so $
\gamma a=a$, because
$\sigma_a\infty =a$, hence $\gamma\in\Gamma_a$, the lemma is proved.

{\bf LEMMA 3.2.} {\it Assume that} $m$ {\it is a compactly supported }
{\it function with bounded variation on}
$[0,\infty )${\it . Let} $a\in A${\it , and let} $Y$ {\it be large enough depending on }
$\Gamma$ {\it and} $m${\it . For} $z,w\in H$ {\it define}
$$M(z,w):=\sum_{\gamma\in\Gamma}m\left({{\left|z-\gamma w\right|^
2}\over {4\hbox{\rm Im$z$$\hbox{\rm Im}\gamma w$}}}\right),\eqno
(3.4)$$
{\it then for} $z\in F_a\left(Y\right)$ {\it we have}
$$M\left(z,z\right)=4\hbox{\rm Im}w\int_0^{\infty}m\left(y^2\right
)dy+O_{\Gamma ,m}\left(1\right),$$
{\it where} $z=\sigma_aw${\it . }

{\it Proof.\/} It follows easily from Lemma 3.1 using (3.2) and
$\sigma_a^{-1}\Gamma_a\sigma_a=B$ that if $Y$ is large enough, then
$$M\left(z,z\right)=\sum_{l=-\infty}^{\infty}m\left({{l^2}\over {
4\hbox{\rm $\hbox{\rm Im}^2w$}}}\right),$$
{\it where} $z=\sigma_aw${\it .\/} Then using the inequality of Koksma (Theorem 5.1 of
[K-N]) we get the lemma.

{\bf LEMMA 3.3.} {\it Let} $u$ {\it be a fixed} $\Gamma${\it -automorphic eigenfunction of the Laplace operator with eigenvalue}
$\lambda =s(s-1)$ {\it satisfying}
$$\int_F\left|u\left(z\right)\right|d\mu_z<\infty ,\eqno (3.5)$$
{\it and let} $\hbox{\rm Re}s={1\over 2}$ {\it or} ${1\over 2}<s\le
1${\it . Denote the}
{\it Fourier expansion of} $u$ {\it by}
$$u\left(\sigma_az\right){\it =}\beta_a\left(0\right)y^s{\it +}\tilde{
\beta}_a\left(0\right)y^{1-s}{\it +}\sum_{n\neq 0}\beta_a\left(n\right
)W_s\left(nz\right){\it .}$$
{\it Introduce the notations}
$$B_u{\it =}\sum_a\beta_a\left(0\right){\it ,}\qquad\tilde {B}_u{\it =}
\sum_a\tilde{\beta}_a\left(0\right){\it .}$$
{\it Assume that} $m$ {\it is a compactly supported function with bounded variation on} $
[0,\infty )$ {\it and}
$$\int_0^{\infty}{{m\left(v\right)}\over {\sqrt v}}dv=0.\eqno (3.
6)$$
{\it Recalling the notation (3.4) we have that}
$$\int_FM(z,z)u(z)d\mu_z=\Sigma_{hyp}+\Sigma_{ell}+\Sigma_{par},\eqno
(3.7)$$
{\it with the definitions}
$$\Sigma_{hyp}:=\sum_{\matrix{_{[\gamma ]}\cr
\gamma\hbox{\rm \ hyperbolic}\cr}
}\left(\int_{C_{\gamma}}udS\right)\int_{-{{\pi}\over 2}}^{{{\pi}\over
2}}m\left({{N\left(\gamma\right)+N\left(\gamma\right)^{-1}-2}\over {
4\cos^2\theta}}\right)f_{\lambda}(\theta ){{d\theta}\over {\cos^2
\theta}},$$
{\it where the summation is over the hyperbolic conjugacy classes of} $
\Gamma${\it ,} $ $$N\left(\gamma\right)$ {\it is}
{\it the norm of (the conjugacy class of)} $\gamma${\it ,} $C_{\gamma}$ {\it is the closed geodesic obtained}
{\it by factorizing the noneuclidean line connecting the fixed points of} $
\gamma$ {\it by the action of the centralizer of} $\gamma$ {\it in} $
\Gamma${\it ,} $dS={{\left|dz\right|}\over y}$ {\it is the hyperbolic arc length},

$$\Sigma_{ell}:=\sum_{p\in P}{{2\pi}\over {m_p}}u\left(p\right)\sum_{
l=1}^{m_p-1}\int_0^{\infty}m\left(\sin^2{{l\pi}\over {m_p}}\sinh^
2r\right)g_{\lambda}(r)\sinh rdr,$$
{\it and for} $s\ne 1$ {\it we write}
$$\Sigma_{par}:=B_u2^{1-s}\zeta\left(1-s\right)\int_0^{\infty}{{m\left
(v\right)}\over {v^{{{1+s}\over 2}}}}dv+\tilde {B}_u2^s\zeta\left
(s\right)\int_0^{\infty}{{m\left(v\right)}\over {v^{{{2-s}\over 2}}}}
dv,$$
{\it for} $s=1$ {\it we write}
$$\Sigma_{par}:=\tilde {B}_u\int_0^{\infty}{{m\left(v\right)}\over {
v^{{1\over 2}}}}\log vdv$$
{\it where} $\zeta$ {\it is the Riemann zeta function. The left-hand side of (3.7) is absolutely convergent and} $
\Sigma_{hyp}$ {\it is a finite sum. }

{\it Proof.\/} This is essentially proved in [B1] in the case $s\ne
1$ and in [I] in the case $s=1$ (since for $s=1$ this follows from the classical Selberg
trace formula), but it is not stated there exactly in this form, so we explain how
it follows from [B1] and from [I].

It follows from Lemma 3.2 and condition (3.6) that $M(z,z)$
is bounded on $H$, hence by (3.5) the left-hand side of (3.7) is
absolutely convergent. Let us write
$$M_{hyp}(z):=\sum_{\matrix{\gamma\in\Gamma\cr
\gamma\hbox{\rm \ hyperbolic}\cr}
}m\left({{\left|z-\gamma z\right|^2}\over {4\hbox{\rm Im$z$$\hbox{\rm Im}
\gamma z$}}}\right),$$
$$M_{ell}(z):=\sum_{\matrix{\gamma\in\Gamma\cr
\gamma\hbox{\rm \ elliptic}\cr}
}m\left({{\left|z-\gamma z\right|^2}\over {4\hbox{\rm Im$z$$\hbox{\rm Im}
\gamma z$}}}\right),$$
$$M_{par}(z):=\sum_{\matrix{\gamma\in\Gamma\cr
\gamma\hbox{\rm \ parabolic}\cr}
}m\left({{\left|z-\gamma z\right|^2}\over {4\hbox{\rm Im$z$$\hbox{\rm Im}
\gamma z$}}}\right).$$
It is clear by Lemma 3.1 and the fact that $\Gamma$ acts
discontinuously on $H$ (see p 40 of [I]) that there are only
finitely many $\gamma\in\Gamma$ for which there is a $z\in H$ such that
the contribution of $\gamma$ to $M_{hyp}(z)$ or $M_{ell}(z)$ is nonzero.
It follows then, on the one hand, that $M_{par}(z)$ is also
bounded on $H$, and on the other hand that

$$\int_FM_{hyp}(z)u(z)d\mu_z=\sum_{\matrix{_{[\gamma ]}\cr
\gamma\hbox{\rm \ hyperbolic}\cr}
}T_{\gamma},\qquad\int_FM_{ell}(z)u(z)d\mu_z=\sum_{\matrix{_{[\gamma
]}\cr
\gamma\hbox{\rm \ elliptic}\cr}
}T_{\gamma}$$
where the summation is over the hyperbolic and elliptic
conjugacy classes of $\Gamma$, respectively, $ $and
$$T_{\gamma}:=\sum_{\delta\in [\gamma ]}\int_Fm\left({{\left|z-\delta
z\right|^2}\over {4\hbox{\rm Im$z$$\hbox{\rm Im}\delta z$}}}\right
)u(z)d\mu_z.$$
Then it follows from (3) and (4) of [B1] (and the
reasoning there is valid also for the case $s=1$) that
$$\int_FM_{hyp}(z)u(z)d\mu_z=\Sigma_{hyp},\qquad\int_FM_{ell}(z)u
(z)d\mu_z=\Sigma_{ell}.$$
Since we have seen that $M_{par}(z)$ is bounded, so
$$\int_FM_{par}(z)u(z)d\mu_z=\Sigma_{par}$$
follows from Lemma 3 of [B1] in the case $s\ne 1$ and from
(10.14) and (10.15) of [I] in the case $s=1$ (since in that
case $u$ is constant), taking into account that
$g_m\left(0\right)=0$ by our condition (3.6) (see (2.2)). The lemma is proved.

{\bf LEMMA 3.4.} {\it Let} $m$ {\it be a compactly supported continuous function on}
$[0,\infty )${\it , Assume that} $h_m\left(r\right)$ {\it (defined in (2.2) and (2.3)) is even, it is holomorphic in the}
{\it strip} $\left|\hbox{\rm Im}r\right|\le{1\over 2}+$$\epsilon$ {\it and} $
h_m\left(r\right)=O\left(\left(1+\left|r\right|\right)^{-2-\epsilon}\right
)$ {\it in this}
{\it strip for some} $\epsilon >0${\it . Then for} $z\in H$ {\it we have (using }
{\it definition (3.4)) that }
$$M(z,z)=\sum_{j=0}^{\infty}h_m\left(t_j\right)\left|u_j\left(z\right
)\right|^2+\sum_a{1\over {4\pi}}\int_{-\infty}^{\infty}h_m\left(r\right
)\left|E_a\left(z,{1\over 2}+ir\right)\right|^2dr,$$
{\it and this expression is absolutely and uniformly convergent on compact subsets of} $
H${\it .}

{\it Proof.\/} This follows from Theorem 7.4 of [I].

\noindent{\bf 4. Proof of the theorem}
\medskip

Let $x$ be a large positive number and let $1<d<{x\over {\log x}}$
(say). We will choose later the parameter $d$ optimally.

Let
$$k\left(y\right)=1\hbox{\rm \ for $0\le y\le x$},\quad k\left(y\right
)=0\hbox{\rm \ for $y>x$},\eqno (4.1)$$
and let $k^{\ast}$ be a smoothed version of $k$, more precisely
$$k^{\ast}\left(y\right):=\int_{-\infty}^{\infty}\eta\left(\tau\right
)k\left(ye^{\tau}\right)d\tau ,\eqno (4.2)$$
where the function $\eta$ will be a smooth even function satisfying
$\eta\left(\tau\right)=0$ for $\left|\tau\right|>d/x$. We define $
\eta$ precisely below.

But before defining $\eta$ let us remark that our goal is to
achieve
$$\int_0^{\infty}\left(k^{\ast}(u)-k(u)\right)u^{-1/2}du=0,\eqno
(4.3)$$
since we would like to apply Lemma 3.3 for this
difference. By (4.2) we have
$$\int_0^{\infty}k^{\ast}(u)u^{-1/2}du=\int_{-\infty}^{\infty}\eta\left
(\tau\right)e^{-\tau /2}d\tau\int_0^{\infty}k\left(v\right)v^{-1/
2}dv,\eqno (4.4)$$
and so we want to have
$$\int_{-\infty}^{\infty}\eta\left(\tau\right)e^{-\tau /2}d\tau =
1.$$
We will take
$$\eta\left(\tau\right):={x\over d}\eta_0\left({x\over d}\tau\right
),\eqno (4.5)$$
where the function $\eta_0$ will be a smooth even function satisfying
$\eta_0\left(\tau\right)=0$ for $\left|\tau\right|>1$. Then
$$\int_{-\infty}^{\infty}\eta\left(\tau\right)e^{-\tau /2}d\tau =
\int_{-\infty}^{\infty}\eta_0\left(\tau\right)e^{-\tau{d\over {2x}}}
d\tau ,\eqno (4.6)$$
so we need
$$\int_{-\infty}^{\infty}\eta_0{}_{}\left(\tau\right)e^{-\tau{d\over {
2x}}}d\tau =1.\eqno (4.7)$$
We now define $\eta$$_0$. First let $\psi_0$ be a given smooth even nonnegative function on the real line such that
$$\psi_0\left(\tau\right)=0\hbox{\rm \ for $\left|\tau\right|>1$}$$
and
$$\int_{-\infty}^{\infty}\psi_0{}_{}\left(\tau\right)d\tau =1.\eqno
(4.8)$$
Then
$$I_{d,x}:=\int_{-\infty}^{\infty}\psi_0{}_{}\left(\tau\right)e^{
-\tau{d\over {2x}}}d\tau =1+O\left(\left({d\over x}\right)^2\right
)\eqno (4.9)$$
with implied absolute constant. If $x$ is large enough, then
clearly $1/2<I_{d,x}<2$. Let
$$\eta_0\left(\tau\right)={{\psi_0{}_{}\left(\tau\right)}\over {I_{
d,x}}}\eqno (4.10)$$
for real $\tau$, then we have (4.7). Note that $\eta_0$ slightly depends on $
d$ and $x$, but we do not
denote it. Formulas (4.9), (4.10), (4.5) define $\eta$, and by
(4.7), (4.4), (4.6) we get (4.3).

Note that by the definitions for any integer $j\ge 0$ we have that
$$\int_{-\infty}^{\infty}\left|\eta^{(j)}\left(\tau\right)\right|
d\tau\ll_j\left({x\over d}\right)^j,\eqno (4.11)$$
and we also have (by (4.8), (4.9), (4.10) and (4.5)) that
$$\int_{-\infty}^{\infty}\eta\left(\tau\right)d\tau =1+O\left(\left
({d\over x}\right)^2\right).\eqno (4.12)$$
So the smoothed version $k^{\ast}$ of $k$ is now defined. As it is
mentioned in Remark 1.2, we will use the decomposition
(1.2), and we will apply Lemma 3.4 for the first term
there, we will apply Lemma 3.3 for the second term. To
apply these lemmas we need estimates for the function
transforms occurring in those lemmas. We give such
estimates in the next three lemmas.

For simplicity introduce the abbreviations $q^{\ast}=q_{k^{\ast}}$, $
g^{\ast}=g_{k^{\ast}}$ and $h^{\ast}=h_{k^{\ast}}$ (see (2.2) and (2.3)).

{\bf LEMMA 4.1.} {\it For every integer} $j\ge 2$ {\it we have for} $
r\ge 1$ {\it that}
$$\left|h^{\ast}\left(r\right)\right|\ll_j{{d^{3/2}}\over x}\left
({x\over {dr}}\right)^j.\eqno (4.13)$$
{\it We also have for every complex} $r$ {\it that}
$$\left|h^{\ast}\left(r\right)\right|\ll x^{{1\over 2}+\left|\hbox{\rm Im$
r$}\right|}\log x.\eqno (4.14)$$
{\it Proof.\/} By (2.2) and (4.1) we have
$$q_k\left(y\right)=2\sqrt {x-y}\hbox{\rm \ for $0\le y\le x$},\quad
q_k\left(y\right)=0\hbox{\rm \ for $y>x$}.\eqno (4.15)$$
It is easy to see by (2.2) and (4.2) that we have
$$q^{\ast}\left(v\right)=\int_{-\infty}^{\infty}\eta\left(\tau\right
)e^{-{{\tau}\over 2}}q_k\left(ve^{\tau}\right)d\tau ,\eqno (4.16)$$
and by the substitution $e^{\mu}=ve^{\tau}$ we can also write
$$q^{\ast}\left(v\right)=\sqrt v\int_{-\infty}^{\infty}\eta\left(
\mu -\log v\right)e^{-{{\mu}\over 2}}q_k\left(e^{\mu}\right)d\mu
.$$
Then for any integer $j\ge 0$ we have on the one hand that
$$\left(q^{\ast}\left(v\right)\right)^{(j)}=\int_{-\infty}^{\infty}
\eta\left(\tau\right)e^{\left(j-{1\over 2}\right)\tau}q^{(j)}_k\left
(ve^{\tau}\right)d\tau ,\eqno (4.17)$$
and we have on the other hand that
$$\left(q^{\ast}\left(v\right)\right)^{(j)}=\sqrt v\int_{-\infty}^{
\infty}{{\sum_{l=0}^jc_{l,j}\eta^{(l)}\left(\mu -\log v\right)}\over {
v^j}}e^{-{{\mu}\over 2}}q_k\left(e^{\mu}\right)d\mu\eqno (4.18)$$
with some constants $c_{l,j}$.

It is clear from (4.15) and (4.16) that for $v\ge xe^{d/x}$ we have
$$q^{\ast}\left(v\right)=0.\eqno (4.19)$$
It is also clear by the same formulas that for $0\le v\le xe^{d/x}$ we have
$$0\le q^{\ast}\left(v\right)\ll\sqrt x.\eqno (4.20)$$
The estimate (4.14) follows at once from (4.19), (4.20), (2.2) and
(2.3).

Assume that $0\le v\le xe^{-2d/x}$. Then for $\eta\left(\tau\right
)\neq 0$ we have
$$ve^{\tau}\le ve^{d/x}\le v+{{x-v}\over 2},$$
since this latter inequality is easily seen to be equivalent to
$$2e^{d/x}-1\le{x\over v},$$
which is true, since ${x\over v}\ge e^{2d/x}$. So for any integer $
j\ge 1$
we have by (4.15) that
$$\left|q^{(j)}_k\left(ve^{\tau}\right)\right|\ll_j\left(x-ve^{\tau}\right
)^{{1\over 2}-j}\ll_j\left(x-v\right)^{{1\over 2}-j},$$
hence by (4.17) we get that
$$\left|\left(q^{\ast}\left(v\right)\right)^{(j)}\right|\ll_j\left
(x-v\right)^{{1\over 2}-j}\eqno (4.21)$$
for $0\le v\le xe^{-2d/x}$ and $j\ge 1$.

Now let $xe^{-2d/x}\le v\le xe^{d/x}$. Then we use (4.18). If the
integrand here is nonzero, then we must have
$\left|\mu -\log v\right|\le d/x$, so
$$xe^{-3d/x}\le e^{\mu}\le xe^{2d/x},$$
hence by (4.15) one has
$$q_k\left(e^{\mu}\right)\ll\sqrt d,$$
and by the upper and lower bounds for $e^{\mu}$ and $v$ one also has
$$x\ll e^{\mu}\ll x,\qquad x\ll v\ll x$$
Using these estimates, by (4.18) and (4.11) we get
for any integer $j\ge 1$ that
$$\left|\left(q^{\ast}\left(v\right)\right)^{(j)}\right|\ll_jd^{{
1\over 2}-j}\eqno (4.22)$$
for $xe^{-2d/x}\le v\le xe^{d/x}$.

We see by (2.2) for every $j\ge 1$ and real $a$ that
$$\left|\left(g^{\ast}\left(a\right)\right)^{(j)}\right|\ll_j\sum_{
l=1}^j\left|\left(q^{\ast}\left(\sinh^2{a\over 2}\right)\right)^{
(l)}\right|e^{l\left|a\right|}.\eqno (4.23)$$
By (2.3) we have by repeated partial integration for every
$j\ge 1$ and $r\ge 1$ that
$$\left|h^{\ast}\left(r\right)\right|\ll_j{1\over {r^j}}\int_{-\infty}^{
\infty}\left|\left(g^{\ast}\left(a\right)\right)^{(j)}\right|da.$$
By (4.19), (4.21), (4.22) and (4.23) we obtain (4.13). The lemma is
proved.

{\bf LEMMA 4.2.} {\it For} ${1\over {100}}\le it\le{1\over 2}$ {\it (say) we have that}
$$h^{\ast}(t)=\sqrt {\pi}{{\Gamma\left(it\right)2^{2it+1}}\over {
\Gamma\left({3\over 2}+it\right)}}x^{{1\over 2}+it}+O\left(x\left
({d\over x}\right)^2\right)+O\left(x^{{1\over 2}}\right).$$
{\it Proof.\/} Assume first $it<{1\over 2}$. By (1.62') of [I] (the function $
F_s(u)$ is
defined by the formulas on p.26., line 7, and (B.23) of [I])
and by [G-R], p 995, 9.113 we have that
$$h^{\ast}(t)={2\over {i\Gamma\left({1\over 2}\pm it\right)}}\int_{
(\sigma )}{{\Gamma\left({1\over 2}\pm it+S\right)\Gamma\left(-S\right
)}\over {\Gamma\left(1+S\right)}}\left(\int_0^{\infty}k^{\ast}(u)
u^Sdu\right)dS,$$
where $it-{1\over 2}<\sigma <0$, so by (4.1) and (4.2) we get
$$h^{\ast}(t)={2\over {i\Gamma\left({1\over 2}\pm it\right)}}\int_{
-\infty}^{\infty}\eta\left(\tau\right)\int_{(\sigma )}{{\Gamma\left
({1\over 2}\pm it+S\right)\Gamma\left(-S\right)}\over {\Gamma\left
(1+S\right)}}{{\left(x/e^{\tau}\right)^{S+1}}\over {S+1}}dSd\tau
.$$
Shifting the integration to the left we get (and
observe that it is also true for $it={1\over 2}$) that
$$h^{\ast}(t)={{4\pi}\over {\Gamma\left({1\over 2}+it\right)}}{{\Gamma\left
(2it\right)}\over {\Gamma\left({3\over 2}+it\right)}}\int_{-\infty}^{
\infty}\eta\left(\tau\right)\left(x/e^{\tau}\right)^{{1\over 2}+i
t}d\tau +O\left(x^{{1\over 2}}\right),$$
and taking into account the properties of $\eta$ and the
duplication formula for the $\Gamma$-function ([I], (B.5)) we
obtain the lemma.

If $m$ is a compactly supported continuous function on
$[0,\infty )$, $\lambda\le 0$ and $T>0$, introduce the notation
$$M_{m,\lambda}\left(T\right):=\int_{-{{\pi}\over 2}}^{{{\pi}\over
2}}m\left({T\over {\cos^2\theta}}\right)f_{\lambda}(\theta ){{d\theta}\over {\cos^
2\theta}}.\eqno (4.24)$$

Note that this kind of transform appeared in [B1],
equation (7) and also in [Hu], equation (31).

We obviously have
$$M_{k^{\ast},\lambda}\left(T\right)=\int_{-\infty}^{\infty}\eta\left
(\tau\right)M_{k,\lambda}\left(Te^{\tau}\right)d\tau .\eqno (4.25)$$
{\bf LEMMA 4.3.} {\it (i) For every} $T>0$ {\it the functions  }
$$M_{k,-{1\over 4}-t^2}\left(T\right),\qquad M_{k^{\ast},-{1\over
4}-t^2}\left(T\right)$$
{\it are entire functions of} $t${\it .}

{\it (ii) Let} $0<\delta <{1\over 2}$ {\it be given. There is a constant }
$A_{\delta}>0$ {\it depending only on} $\delta$ {\it such that for every} $
t$ {\it satisfying} ${1\over 4}+t^2\ge 0$ {\it or} $\left|\hbox{\rm Re}\left
(it\right)\right|\le{1\over 2}-\delta$
{\it one has the following statements with the notation }
$\lambda =-{1\over 4}-t^2$:

{\it For} $T\ge xe^{d/x}$ {\it we ha\/}v{\it e that}
$$M_{k^{\ast},\lambda}\left(T\right)=M_{k,\lambda}\left(T\right)=
0,\eqno (4.26)$$
{\it for} $xe^{-2d/x}\le T\le xe^{d/x}$ {\it we have that}
$$\left|M_{k^{\ast},\lambda}\left(T\right)\right|+\left|M_{k,\lambda}\left
(T\right)\right|\ll_{\delta}\left(1+\left|t\right|\right)^{A_{\delta}}\left
({d\over x}\right)^{1/2},\eqno (4.27)$$
{\it and for} $0<T\le xe^{-2d/x}$  {\it we have that}
$$M_{k^{\ast},\lambda}\left(T\right)-M_{k,\lambda}\left(T\right)\ll_{
\delta}{{\left(1+\left|t\right|\right)^{A_{\delta}}d^2}\over {T^{
1/2}\left(x-T\right)^{3/2}}}.$$
{\it Proof.\/} Note first that we can give an explicit formula for $
f_{\lambda}$, namely
$$f_{\lambda}(\theta )=F\left(\matrix{{1\over 4}+{{it}\over 2},{1\over
4}-{{it}\over 2}\cr
{1\over 2}\cr}
;-{{\sin^2\theta}\over {\cos^2\theta}}\right)\eqno (4.28)$$
for $\theta\in (-{{\pi}\over 2},{{\pi}\over 2})$, where $\lambda
=-{1\over 4}-t^2$. This can be proved in
the following way. One has
$${1\over {\pi^{{1\over 2}}}}\int_{-{{\pi}\over 2}}^{{{\pi}\over
2}}F\left(\matrix{{1\over 4}+{{it}\over 2},{1\over 4}-{{it}\over
2}\cr
{1\over 2}\cr}
;-{{\sin^2\theta}\over {\cos^2\theta}}\right)\cos^{2s}\theta{{d\theta}\over {\cos^
2\theta}}={{\Gamma (s-{1\over 4}+{{it}\over 2})\Gamma (s-{1\over
4}-{{it}\over 2})}\over {\Gamma^2(s)}}$$
for $\hbox{\rm Re$s>{1\over 2}$}$, as one can see by the substitution $
y={{\sin^2\theta}\over {\cos^2\theta}}$
and by [G-R], p 807, 7.512.10. By Lemma 11 of [B1] it
follows that
$$\int^{{{\pi}\over 2}}_0\left(F\left(\matrix{{1\over 4}+{{it}\over
2},{1\over 4}-{{it}\over 2}\cr
{1\over 2}\cr}
;-{{\sin^2\theta}\over {\cos^2\theta}}\right)-f_{\lambda}(\theta
)\right)\cos^n\theta d\theta =0$$
for every nonnegative integer $n$, which easily implies
(4.28).

We can see part (i) at once.

To show part (ii) note that using (4.24), (4.28) and the substitution
$$Y=\log{1\over {\cos^2\theta}}$$
one has for $T\le x$ that

$$M_{k,\lambda}\left(T\right)=\int_0^{\log x-\log T}F\left(\matrix{{
1\over 4}+{{it}\over 2},{1\over 4}-{{it}\over 2}\cr
{1\over 2}\cr}
;1-e^Y\right){{e^YdY}\over {\sqrt {e^Y-1}}},.\eqno (4.29)$$
and for $T>x$ we have $M_{k,\lambda}\left(T\right)=0$. Then (4.26) is obvious by
(4.25). We also have by (4.29) and Lemma 5.1 for $T\le x$ that
$$M_{k,\lambda}\left(T\right)\ll_{\delta}\left(1+\left|t\right|\right
)^{A_{\delta}}\left({x\over T}-1\right)^{1/2},\eqno (4.30)$$
and then (using also (4.25)) (4.27) follows.

By (4.25), (4.12) and since $\eta$ is even, we have that
$$M_{k^{\ast},\lambda}\left(T\right)-M_{k,\lambda}\left(T\right)\eqno
(4.31)$$
equals
$$\int_0^{\infty}\eta\left(\tau\right)\left(M_{k,\lambda}\left(Te^{
\tau}\right)+M_{k,\lambda}\left(Te^{-\tau}\right)-2M_{k,\lambda}\left
(T\right)\right)d\tau +O\left(\left({d\over x}\right)^2\right)M_{
k,\lambda}\left(T\right).\eqno (4.32)$$

Assuming $Te^{\tau}\le x$ by (4.29) we have that
$$M_{k,\lambda}\left(Te^{\tau}\right)+M_{k,\lambda}\left(Te^{-\tau}\right
)-2M_{k,\lambda}\left(T\right)=$$
$$=\int_{\log{x\over T}}^{\log{x\over T}+\tau}\left({{F\left(\matrix{{
1\over 4}+{{it}\over 2},{1\over 4}-{{it}\over 2}\cr
{1\over 2}\cr}
;1-e^Y\right)e^Y}\over {\sqrt {e^Y-1}}}-{{F\left(\matrix{{1\over
4}+{{it}\over 2},{1\over 4}-{{it}\over 2}\cr
{1\over 2}\cr}
;1-e^{Y-\tau}\right)e^{Y-\tau}}\over {\sqrt {e^{Y-\tau}-1}}}\right
)dY.$$
For $0<T\le xe^{-2d/x}$ and $\left|\tau\right|\le d/x$ we then have by the
mean-value theorem and by Lemma 5.1 that
$$M_{k,\lambda}\left(Te^{\tau}\right)+M_{k,\lambda}\left(Te^{-\tau}\right
)-2M_{k,\lambda}\left(T\right)\ll_{\delta}\left(1+\left|t\right|\right
)^{A_{\delta}}\tau^2\left({x\over T}\right)^{1/2}\left(1+{1\over {\log{
x\over T}}}\right)^{3/2}.$$
So by (4.30), (4.11) with $j=0$ and by (4.31), (4.32), using that $\eta\left
(\tau\right)=0$ for $\left|\tau\right|>d/x$
we get for $0<T\le xe^{-2d/x}$ that
$$M_{k^{\ast},\lambda}\left(T\right)-M_{k,\lambda}\left(T\right)\ll_{
\delta}\left(1+\left|t\right|\right)^{A_{\delta}}\left({d\over x}\right
)^2\left({x\over T}\right)^{1/2}\left({x\over {x-T}}\right)^{3/2}
.$$
The lemma is proved.

{\bf REMARK 4.4.} We now compare (4.28) to formulas (86) and (87) of [F]. Indeed, by the notation used
there we see that
$$\hbox{\rm $N_{{1\over 2}+it,0}^0\left(\theta +{{\pi}\over 2}\right
)$ and  $N_{{1\over 2}-it,0}^0\left(\theta +{{\pi}\over 2}\right)$}$$
are solutions of our differential equation (2.1) for
$-{{\pi}\over 2}<\theta <0$. Using (87) of [F] and the quadratic
transformation [G-R], p 999, 9.134.1 one sees that
$$N_{{1\over 2}+it,0}^0\left(\theta +{{\pi}\over 2}\right)=\left({{\cos^
2\theta}\over {\sin^2\theta}}\right)^{{1\over 4}+{{it}\over 2}}F\left
(\matrix{{1\over 4}+{{it}\over 2},{3\over 4}+{{it}\over 2}\cr
1+it\cr}
;-{{\cos^2\theta}\over {\sin^2\theta}}\right)$$
for $-{{\pi}\over 2}<\theta <0$. Using the same relation with $-t$ in place
of $t$ we see by [G-R], p 999, 9.132.2 that the right-hand
side of (4.28) is a linear combination of $N_{{1\over 2}+it,0}^0\left
(\theta +{{\pi}\over 2}\right)$ and
$N_{{1\over 2}-it,0}^0\left(\theta +{{\pi}\over 2}\right)$ on the interval $
-{{\pi}\over 2}<\theta <0$, hence it is also a solution of
(2.1). Since it is an even function, it gives a solution of (2.1)
on the whole interval $-{{\pi}\over 2}<\theta <{{\pi}\over 2}$.

Continuing the proof of Theorem 1.1 let
$$m_1\left(v\right)=k^{\ast}\left(v\right),\qquad m_2\left(v\right
)=k\left(v\right)-k^{\ast}\left(v\right).\eqno (4.33)$$
For $m_1\left(v\right)$ we will apply Lemma 3.4, and for $m_2\left(
v\right)$ we
will apply Lemma 3.3.

We can see e.g. by (1.62') of [I] (the function $F_s(u)$ is
defined by the formulas on p.26., line 7, and (B.23) of [I])
and by Lemma 6.2 of [B2] that the conditions of Lemma 3.4 are satisfied writing $
m_1$
in place of $m$, hence for $z\in H$ we have that
$$\sum_{\gamma\in\Gamma}m_1\left({{\left|z-\gamma z\right|^2}\over {
4\hbox{\rm Im$z$$\hbox{\rm Im}\gamma z$}}}\right)\eqno (4.34)$$
equals
$$\sum_{j=0}^{\infty}h_{m_1}\left(t_j\right)\left|u_j\left(z\right
)\right|^2+\sum_a{1\over {4\pi}}\int_{-\infty}^{\infty}h_{m_1}\left
(r\right)\left|E_a\left(z,{1\over 2}+ir\right)\right|^2dr.$$
Then applying Lemma 4.1 (we apply (4.13) with $j=2$ for
$1\le r<{x\over d}$, we apply (4.13) with $j=3$ for $r\ge{x\over
d}$, finally we
apply (4.14) for $\left|\hbox{\rm Re$r$}\right|<1$, $\left|\hbox{\rm Im$
r$}\right|<{1\over {100}}$), Lemma 4.2
and [I], Proposition 7.2, for every $z\in H$
satisfying $f\left(z\right)\neq 0$ and for every $\epsilon >0$ we get that (4.34) equals
$$\sum_{j,\>it_j>0}\sqrt {\pi}{{\Gamma\left(it_j\right)2^{2it_j+1}}\over {
\Gamma\left({3\over 2}+it_j\right)}}x^{{1\over 2}+it_j}\left|u_j\left
(z\right)\right|^2+O_{f,\epsilon}\left(x^{\epsilon}\left({x\over {\sqrt
d}}+x^{{1\over 2}+{1\over {100}}}+{{d^2}\over x}\right)\right),\eqno
(4.35)$$
where we took into account that $f$ is compactly supported on $F$.

Let $f\left(z\right)$ be as in the theorem, and consider the integral
$$\int_Ff\left(z\right)\left(\sum_{\gamma\in\Gamma}m_2\left({{\left
|z-\gamma z\right|^2}\over {4\hbox{\rm Im$z$$\hbox{\rm Im}\gamma
z$}}}\right)\right)d\mu_z.\eqno (4.36)$$
By (4.3) and Lemma 3.2 we see that the function in the
bracket here is bounded. Then it follows from the
Spectral Theorem ([I], Theorems 4.7 and 7.3) that (4.36)
equals
$$\sum_{j=0}^{\infty}\left(f,u_j\right)\int_Fu_j\left(z\right)\left
(\sum_{\gamma\in\Gamma}m_2\left({{\left|z-\gamma z\right|^2}\over {
4\hbox{\rm Im$z$$\hbox{\rm Im}\gamma z$}}}\right)\right)d\mu_z+\eqno
(4.37)$$
$$+\sum_a{1\over {4\pi}}\int_{-\infty}^{\infty}\left(f,E_a\left(\ast
,{1\over 2}+ir\right)\right)\int_FE_a\left(z,{1\over 2}+ir\right)\left
(\sum_{\gamma\in\Gamma}m_2\left({{\left|z-\gamma z\right|^2}\over {
4\hbox{\rm Im$z$$\hbox{\rm Im}\gamma z$}}}\right)\ \right)d\mu_zd
r.$$
By (4.3) we see that the conditions of Lemma 3.3 are
satisfied writing $m_2$ in place of $m$ and writing $u=u_j$ or
$u=E_a\left(\ast ,{1\over 2}+ir\right)$. By Lemma 6.3, Lemma 5.2, (4.2), (4.11)
with $j=0$, using also (6.28) of [I] and a convexity bound for the Riemann zeta function we get
that after applying Lemma 3.3 the contribution of $\Sigma_{ell}$ and
$\Sigma_{par}$ to (4.37) is $O_{f,\epsilon}\left(x^{{1\over 2}+\epsilon}\right
)$ for every $\epsilon >0$.

Therefore, for a hyperbolic $\gamma\in\Gamma$ introducing the notation
$$T\left(\gamma\right)={{N\left(\gamma\right)+N\left(\gamma\right
)^{-1}-2}\over 4}$$
and recalling (4.24) we get that (4.36) equals
$$O_{f,\epsilon}\left(x^{{1\over 2}+\epsilon}\right)+\sum_{\matrix{_{
[\gamma ]}\cr
\gamma\hbox{\rm \ hyperbolic}\cr}
}\left(\Sigma_1\left(\gamma\right)+\Sigma_2\left(\gamma\right)\right
)\eqno (4.38)$$
with the notations
$$\Sigma_1\left(\gamma\right):=\sum_{j=0}^{\infty}\left(f,u_j\right
)\left(\int_{C_{\gamma}}u_j^{Y_{\Gamma}}dS\right)M_{m_2,-{1\over
4}-t_j^2}\left(T\left(\gamma\right)\right)+$$
$$+\sum_a{1\over {4\pi}}\int_{-\infty}^{\infty}\left(f,E_a\left(\ast
,{1\over 2}+ir\right)\right)\left(\int_{C_{\gamma}}E_a^{Y_{\Gamma}}\left
(\ast ,{1\over 2}+ir\right)dS\right)M_{m_2,-{1\over 4}-r^2}\left(
T\left(\gamma\right)\right)dr,\eqno (4.39)$$
$$\Sigma_2\left(\gamma\right):=\sum_{j=0}^{\infty}\left(f,u_j\right
)\left(\int_{C_{\gamma}}\left(u_j-u_j^{Y_{\Gamma}}\right)dS\right
)M_{m_2,-{1\over 4}-t_j^2}\left(T\left(\gamma\right)\right)+$$
$$+\sum_a{1\over {4\pi}}\int_{-\infty}^{\infty}\left(f,E_a\left(\ast
,{1\over 2}+ir\right)\right)\left(\int_{C_{\gamma}}D_a^{Y_{\Gamma}}\left
(\ast ,{1\over 2}+ir\right)dS\right)M_{m_2,-{1\over 4}-r^2}\left(
T\left(\gamma\right)\right)dr,$$
where for simplicity we wrote
$$D_a^{Y_{\Gamma}}\left(\ast ,{1\over 2}+ir\right):=E_a\left(\ast
,{1\over 2}+ir\right)-E_a^{Y_{\Gamma}}\left(\ast ,{1\over 2}+ir\right
).$$
Observe that for any $a\in A$ and any $y>0$ we have by [I],
(6.22) and (6.27) that
$$\sum_c{1\over {4\pi}}\int_{-\infty}^{\infty}\left(f,E_c\left(\ast
,{1\over 2}+ir\right)\right)\left(\delta_{ca}y^{{1\over 2}+ir}+\phi_{
c,a}\left({1\over 2}+ir\right)y^{{1\over 2}-ir}\right)M_{m_2,-{1\over
4}-r^2}\left(T\left(\gamma\right)\right)dr$$
equals
$${1\over {2\pi}}\int_{-\infty}^{\infty}\left(f,E_a\left(\ast ,{1\over
2}-ir\right)\right)y^{{1\over 2}-ir}M_{m_2,-{1\over 4}-r^2}\left(
T\left(\gamma\right)\right)dr.$$
We then see using the notations of Lemma 6.2 (taking into account also that $
E_a\left(\ast ,{1\over 2}-ir\right)$ and
$E_a\left(\ast ,{1\over 2}+ir\right)$ are conjugates of each other for real $
r$) that
$$\Sigma_2\left(\gamma\right)=\sum_{j=0}^{\infty}\left(f,u_j\right
)\left(\int_{C_{\gamma}}\left(u_j-u_j^{Y_{\Gamma}}\right)dS\right
)M_{m_2,-{1\over 4}-t_j^2}\left(T\left(\gamma\right)\right)+$$
$$+\sum_a{1\over {2\pi i}}\int_{\left({1\over 2}\right)}\left(\int_
Ff\left(z\right)E_a\left(z,s\right)d\mu_z\right)\left(\int_{C_{\gamma}}
A_a\left(\ast ,s\right)dS\right)M_{m_2,s(s-1)}\left(T\left(\gamma\right
)\right)ds.\eqno $$
Since $M_{m_2,s(s-1)}\left(T\left(\gamma\right)\right)$ is analytic in $
s$ by Lemma 4.3, so
shifting the line of integration to the right, using
Lemma 6.4 to see that the residues cancel out and using
also Lemma 6.3 (iii) and Lemma 4.3 we get that
$$\Sigma_2\left(\gamma\right)=\left(f,u_0\right)\left(\int_{C_{\gamma}}\left
(u_0-u_0^{Y_{\Gamma}}\right)dS\right)M_{m_2,0}\left(T\left(\gamma\right
)\right)+$$
$$+\sum_a{1\over {2\pi i}}\int_{\left(1-\delta\right)}\left(\int_
Ff\left(z\right)E_a\left(z,s\right)d\mu_z\right)\left(\int_{C_{\gamma}}
A_a\left(\ast ,s\right)dS\right)M_{m_2,s(s-1)}\left(T\left(\gamma\right
)\right)ds,\eqno (4.40)$$
where $\delta >0$ is a small number chosen in such a way
that $1-\delta >S_l$ for every $l\in L$ satisfying $S_l<1$.

Choosing $\delta$ small enough in terms of $\epsilon$, applying Lemma
4.3, Lemma 6.3, Lemma 6.2, using (4.38), (4.39) and (4.40) we
get that (4.36) equals
$$O_{f,\epsilon}\left(x^{{1\over 2}+\epsilon}\right)+$$
$$O_{f,\epsilon}\left(x^{\epsilon}\sum_{\matrix{_{[\gamma ]}\cr
\gamma\hbox{\rm \ hyperbolic},\,T\left(\gamma\right)\le xe^{-2d/x}\cr}
}{{d^2\log N\left(\gamma\right)}\over {T\left(\gamma\right)^{1/2}\left
(x-T\left(\gamma\right)\right)^{3/2}}}\right)+$$
$$O_{f,\epsilon}\left(x^{\epsilon}\left({d\over x}\right)^{1/2}\sum_{\matrix{_{
[\gamma ]}\cr
\gamma\hbox{\rm \ hyperbolic},\,xe^{-2d/x}\le T\left(\gamma\right
)\le xe^{d/x}\cr}
}\log N\left(\gamma\right)\right).$$
From the prime geodesic theorem (Theorem 10.5 of [I]) we then get assuming
$$d\ge x^{3/4}$$
that (4.36) equals
$$O_{f,\epsilon}\left(\left(x^{\epsilon}\right)\left(x^{{1\over 2}}
+{{d^2}\over x}+{{d^{3/2}}\over {\sqrt x}}\right)\right).\eqno (4
.41)$$
Then by (4.33), (4.34), (4.36), (4.41) and (4.35) we get
for $x^{3/4}\le d\le{x\over {\log x}}$ that
$$\int_Ff\left(z\right)\left(\sum_{\gamma\in\Gamma}k\left({{\left
|z-\gamma z\right|^2}\over {4\hbox{\rm Im$z$$\hbox{\rm Im}\gamma
z$}}}\right)\right)d\mu_z$$
equals
$$\int_Ff\left(z\right)\left(\sum_{j,\>it_j>0}\sqrt {\pi}{{\Gamma\left
(it_j\right)2^{2it_j+1}}\over {\Gamma\left({3\over 2}+it_j\right)}}
x^{{1\over 2}+it_j}\left|u_j\left(z\right)\right|^2\right)d\mu_z+
O_{f,\epsilon}\left(\left(x^{\epsilon}\right)\left({x\over {\sqrt
d}}+{{d^{3/2}}\over {\sqrt x}}\right)\right).$$
Choosing $d=x^{3/4}$ and $x={{X-2}\over 4}$ we get the theorem.

\noindent$ ${\bf 5. Lemmas on special functions}
\medskip

{\bf LEMMA 5.1.} {\it Let} $0<\delta <{1\over 2}$ {\it be given. There is a constant }
$A_{\delta}>0$ {\it depending only on} $\delta$ {\it such that for every} $
t$ {\it satisfying} ${1\over 4}+t^2\ge 0$ {\it or} $\left|\hbox{\rm Re}\left
(it\right)\right|\le{1\over 2}-\delta$ {\it and for every} $X\ge
0$
{\it one has that}
$$\left|F\left(\matrix{{1\over 4}+{{it}\over 2},{1\over 4}-{{it}\over
2}\cr
{1\over 2}\cr}
;-X\right)\right|+\left|\left(1+X\right){d\over {dX}}F\left(\matrix{{
1\over 4}+{{it}\over 2},{1\over 4}-{{it}\over 2}\cr
{1\over 2}\cr}
;-X\right)\right|\le A_{\delta}\left(1+\left|t\right|\right)^{A_{
\delta}}.$$
{\it Proof.\/} Note that
$${d\over {dX}}F\left(\matrix{{1\over 4}+{{it}\over 2},{1\over 4}
-{{it}\over 2}\cr
{1\over 2}\cr}
;-X\right)=-2\left({1\over {16}}+{{t^2}\over 4}\right)F\left(\matrix{{
5\over 4}+{{it}\over 2},{5\over 4}-{{it}\over 2}\cr
{3\over 2}\cr}
;-X\right)$$
and so
$$\left(1+X\right){d\over {dX}}F\left(\matrix{{1\over 4}+{{it}\over
2},{1\over 4}-{{it}\over 2}\cr
{1\over 2}\cr}
;-X\right)=-2\left({1\over {16}}+{{t^2}\over 4}\right)F\left(\matrix{{
1\over 4}+{{it}\over 2},{1\over 4}-{{it}\over 2}\cr
{3\over 2}\cr}
;-X\right),$$
where we used the third line of [G-R], p 998, 9.131.1.
Then it is trivial by [G-R], p 995, 9.111 that the
statement of the lemma is true for $\left|t\right|<1/10$ (say). The
satement of the lemma is also trivial for every $t$ and
for $\left|X\right|<{1\over {10\left(1+\left|t\right|\right)^2}}$ (say) by estimating trivially the series definition of the
hypergeometric function.

For $j=0,1$ and $X>0$ one has that
$$F\left(\matrix{{1\over 4}+j+{{it}\over 2},{1\over 4}+j-{{it}\over
2}\cr
{1\over 2}+j\cr}
;-X\right)$$
equals the sum of
$${{\Gamma\left({1\over 2}+j\right)\Gamma\left(it\right)\Gamma\left
(1-it\right)}\over {\Gamma\left({1\over 4}\pm{{it}\over 2}\right)
\Gamma\left({3\over 4}-{{it}\over 2}\right)\Gamma\left({1\over 4}
+j+{{it}\over 2}\right)}}\int_0^1y^{-{1\over 4}-{{it}\over 2}}\left
(1-y\right)^{-{3\over 4}-{{it}\over 2}}\left(X+y\right)^{-{1\over
4}-j+{{it}\over 2}}dy$$
and the same expression writing $-t$ in place of $t$, this
follows from [G-R], p 999, 9.132.2 and [G-R], p 995, 9.111.
Then the statement follows for the case $\left|t\right|\ge 1/10$,
$\left|X\right|\ge{1\over {10\left(1+\left|t\right|\right)^2}}$. The lemma is proved.

{\bf LEMMA 5.2.} {\it Let $x>0$ and let the function $k$ be defined
by (4.1). Let $0<a<1$. There is an absolute constant} $A_0>0$ {\it such }
{\it that for every} $t$ {\it satisfying} ${1\over 4}+t^2\ge 0$ {\it one has, writing }
$\lambda =-{1\over 4}-t^2$ {\it that}
$$\left|\int_0^{\infty}k\left(a\sinh^2r\right)g_{\lambda}(r)\sinh
rdr\right|\le A_0\left(1+\left|t\right|\right)^{A_0}\left(1+{x\over
a}\right)^{1/2}.$$
{\it Proof.\/} Note first that we can give an explicit formula for $
g_{\lambda}$, namely
$$g_{\lambda}(r)=F\left(\matrix{{3\over 4}+{{it}\over 2},{3\over
4}-{{it}\over 2}\cr
1\cr}
;-\sinh^2r\right)\cosh r\eqno (5.1)$$
for $r\ge 0$, where $\lambda =-{1\over 4}-t^2$. This can be proved in
the following way. One has
$$\int_0^{\infty}F\left(\matrix{{3\over 4}+{{it}\over 2},{3\over
4}-{{it}\over 2}\cr
1\cr}
;-\sinh^2r\right)\cosh r\sinh^{1-2s}rdr={{\Gamma (s-{1\over 4}\pm{{
it}\over 2})\Gamma (1-s)}\over {2\Gamma ({3\over 4}\pm{{it}\over
2})\Gamma (s)}}$$
for ${1\over 2}<\hbox{\rm Re$s<1$}$, as one can see by the substitution $
y=\sinh^2r$
and by [G-R], p 806, 7.511. By Lemma 11 of [B1] it
follows that
$$\int_0^{\infty}\left(g_{\lambda}(r)-F\left(\matrix{{3\over 4}+{{
it}\over 2},{3\over 4}-{{it}\over 2}\cr
1\cr}
;-\sinh^2r\right)\cosh r\right)\sinh^{1-2s}rdr=0$$
for every ${1\over 2}<\hbox{\rm Re$s<1$}$, which easily implies
(5.1).

Then by the substitution $y=\sinh^2r$ and by the
particular shape of $k$ we see that
$$\int_0^{\infty}k\left(a\sinh^2r\right)g_{\lambda}(r)\sinh rdr$$
equals
$${1\over 2}\int_0^{x/a}F\left(\matrix{{3\over 4}+{{it}\over 2},{
3\over 4}-{{it}\over 2}\cr
1\cr}
;-y\right)dy.$$
Since the integrand here equals
$$\left(1+y\right)^{-1/2}{1\over {\pi}}\int_0^1q^{-1/2}\left(1-q\right
)^{-1/2}F\left(\matrix{{1\over 4}+{{it}\over 2},{1\over 4}-{{it}\over
2}\cr
1/2\cr}
;-qy\right)dq$$
by [G-R], p 807, 7.512.11 and the third line of [G-R], p 998,
9.131.1, so using Lemma 5.1 the present lemma is proved.

\noindent{\bf 6. Lemmas on automorphic functions}
\medskip

{\bf LEMMA 6.1.} {\it There is a constant} $A_{\Gamma}>0$ {\it depending only on} $
\Gamma$ {\it such that if} $\gamma\in\Gamma$ {\it is hyperbolic and} $
z\in H$ {\it is a point on the noneuclidean line connecting the fixed }
{\it points of} $\gamma${\it , then for every} $a\in A$ {\it one has  }
$$\sigma_a^{-1}\delta z\notin P\left(A_{\Gamma}\sqrt {N\left(\gamma\right
)}\right)$$
{\it for every} $\delta\in\Gamma${\it .}

{\it Proof.\/} By the conditions on $\gamma$ and $z$ we have (see [I], p 19) that
$$\rho (z,\gamma z)=\log N\left(\gamma\right),$$
where $\rho$ is the distance function on $H$ defined in (1.2) of
[I]. Then by (1.3) of [I] we have that
$$u(z,\gamma z)={{N\left(\gamma\right)+N\left(\gamma\right)^{-1}-
2}\over 4},$$
and so for every $\delta\in\Gamma$ one has
$$u\left(\delta z,\left(\delta\gamma\delta^{-1}\right)\delta z\right
)={{N\left(\gamma\right)+N\left(\gamma\right)^{-1}-2}\over 4}.\eqno
(6.1)$$
Now, $\delta\gamma\delta^{-1}$ is hyperbolic, so it is not an element of $
\Gamma_a$.
Assume that $A_{\Gamma}$ is large enough and $\sigma_a^{-1}\delta
z\in P\left(A_{\Gamma}\sqrt {N\left(\gamma\right)}\right)$,
then by Lemma 3.1 we have
$$u\left(\delta z,\left(\delta\gamma\delta^{-1}\right)\delta z\right
)\ge D_{\Gamma}A_{\Gamma}^2N\left(\gamma\right).$$
If $A_{\Gamma}$ is large enough, then this contradicts (6.1), the
lemma is proved.

{\bf LEMMA 6.2.} {\it For} $a\in A$ {\it and any complex} $s$ {\it write}
$$A_a\left(z,s\right){\it :}{\it =}\left(\hbox{\rm Im}\sigma_a^{-
1}z\right)^{1-s}\hbox{\rm \ for }z\in F_a\left(Y_{\Gamma}\right){\it ,}$$
$$A_a\left(z,s\right):=0\hbox{\rm \ for $ $}z\in F\setminus F_a\left
(Y_{\Gamma}\right),$$
{\it finally let} $A_a\left(\gamma z,s\right)=A_a\left(z,s\right)$ {\it for} $
\gamma\in\Gamma$ {\it and} $z\in F${\it . Let} $\gamma\in\Gamma$
{\it be hyperbolic, then for} $0\le\sigma :=\hbox{\rm Re$s\le 1$}$ {\it one has (using the notations of Lemma 3.3) that}
$$\int_{C_{\gamma}}\left(1+\sum_{a\in A}\left|A_a\left(\ast ,s\right
)\right|\right)dS\ll\left(N\left(\gamma\right)\right)^{{{1-\sigma}\over
2}}\log N\left(\gamma\right).$$
{\it Proof.\/} It is well-known that $\int_{C_{\gamma}}1dS\le\log
N\left(\gamma\right)$, so it is
enough to show that
$$\left|A_a\left(z,s\right)\right|\ll\left(N\left(\gamma\right)\right
)^{{{1-\sigma}\over 2}}$$
for every $a\in A$ and every $z\in H$ lying on the noneuclidean line connecting the fixed
points of $\gamma$, and this follows from Lemma 6.1. The
lemma is proved.

{\bf LEMMA 6.3.} {\it Let the function} $f$ {\it be as in Theorem 1.1.}

{\it (i) For integers} $j\ge 0$ {\it one has that}
$$\sup_{z\in F}\left|u_j^{Y_{\Gamma}}\left(z\right)\right|\ll\left
(1+\left|t_j\right|\right)^C$$
{\it and for any} $a\in A$ {\it and} $R>0$ {\it one has that}
$$\sup_{z\in F}\int_{-R}^R\left|E_a^{Y_{\Gamma}}\left(z,{1\over 2}
+ir\right)\right|^2dR\ll\left(1+R\right)^C$$
{\it with some absolute constant} $C${\it .}

{\it (ii) For every positive integer} $K$ {\it one has for} $j\ge
0$ {\it that}
$$\left|\left(f,u_j\right)\right|\ll_{f,K}\left(1+\left|t_j\right
|\right)^{-K}.$$
{\it (iii) For any} $a\in A$ {\it the function}
$$\int_Ff\left(z\right)E_a\left(z,s\right)d\mu_z$$
{\it is meromorphic for} ${1\over 2}\le\hbox{\rm Re$s$}\le 2$ {\it having poles only at the }
{\it points} $\left\{S_l:\;l\in L\right\}${\it , and for every positive integer} $
K$ {\it and every} ${1\over 2}\le\sigma\le 2$ {\it one has  }
$$\int_{-\infty}^{-1}\left|\int_Ff\left(z\right)E_a\left(z,\sigma
+ir\right)d\mu_z\right|^2r^{2K}dr\ll_{f,K}1$$
{\it and}
$$\int_1^{\infty}\left|\int_Ff\left(z\right)E_a\left(z,\sigma +ir\right
)d\mu_z\right|^2r^{2K}dr\ll_{f,K}1.$$
{\it Proof.\/} Part (i) follows e.g. from [I], Proposition 7.2 and
formulas (9.13), (8.1), (8.2), (8.5), (8.6). For parts (ii) and (iii) we
use that the Laplace operator is self-adjoint (see (4.2) of
[I]), and we apply repeatededly (4.2) of [I]. Part (ii) follows
at once in this way. Part (iii) also follows in this way if we can show that there is an absolute constant $
K_0>0$
such that
$$\int_{-\infty}^{-1}{{\int_F\left|E_a^Y\left(z,\sigma +ir\right)\right
|^2d\mu_z}\over {r^{2K_0}}}dr+\int^{\infty}_1{{\int_F\left|E_a^Y\left
(z,\sigma +ir\right)\right|^2d\mu_z}\over {r^{2K_0}}}dr\ll_Y1\eqno
(6.2)$$
for ${1\over 2}<\sigma\le 2$ and $Y\ge 1${\it .\/} (We can assume indeed $
\sigma >{1\over 2}$,
since the case $\sigma ={1\over 2}$ of the last statement of (iii)
follows from the ${1\over 2}<\sigma\le 2$ case of that statement by continuity, since the upper bound is uniform in $
\sigma$.)

Statement (6.2) can be
deduced from the Maass-Selberg relations in the following
way. One first shows by (6.31) of [I] (using the notations
of that book) that
$$\left|\phi_{a,a}\left(\sigma +ir\right)\right|\ll 1$$
for ${1\over 2}<\sigma\le 2$, $\left|r\right|\ge 1$ and for any cusp $
a$, then the same
estimate follows for $\phi_{a,b}\left(\sigma +ir\right)$ for any two cusps $
a$, $b$.
Finally, still by (6.31) of [I], we get
$$\sum_b\left|\phi_{a,b}\left(\sigma +ir\right)\right|^2\le 1+O\left
(\sigma -{1\over 2}\right)$$
for ${1\over 2}<\sigma\le 2$, $\left|r\right|\ge 1$ and for any cusp $
a${\it .\/} By the Hadamard
inequality (see e.g. Corollary 7.8.2 of [H-J]) we then see that for
the determinant $\phi =\det$$\left(\phi_{a,b}\right)$ we have that
$$\left|\phi\left(\sigma +ir\right)\right|^2\le\prod_a\left(\sum_
b\left|\phi_{a,b}\left(\sigma +ir\right)\right|^2\right)\le\left(
1+O\left(\sigma -{1\over 2}\right)\right)\sum_b\left|\phi_{a_0,b}\left
(\sigma +ir\right)\right|^2$$
for any fixed cusp $a_0$, and combining it with Propositions
12.7 and 12.8 of [He] we get that
$$\sum_b\left|\phi_{a_0,b}\left(\sigma +ir\right)\right|^2\ge 1-O\left
(\left(\sigma -{1\over 2}\right)\omega\left(r\right)\right)$$
for ${1\over 2}<\sigma\le 2$, $\left|r\right|\ge 1$ and for any cusp $
a$$_0$ with $\omega$ defined
in Proposition 12.7 of [He]{\it .\/} Then using again (6.31) of [I]
and Proposition 12.7 of [He] we get (6.2). The lemma is proved.

{\bf LEMMA 6.4.} {\it For any} $a\in A$ {\it and} $l\in L$ {\it such that} ${
1\over 2}<S_l<1$ {\it we have that}
$$\sum_{j\ge 0,\,s_j=S_l}\left(f,u_j\right)\beta_{a,j}\left(0\right
)=\int_Ff\left(z\right)\hbox{\rm res}_{s=S_l}E_a\left(z,s\right)d
\mu_z.\eqno (6.3)$$
{\bf REMARK 6.5.} Note that since $f$ could be any function
satisfying the conditions of Theorem 1.1, so we could
easily get
$$\sum_{j\ge 0,\,s_j=S_l}\overline {u_j\left(z\right)}\beta_{a,j}\left
(0\right)=\hbox{\rm res}_{s=S_l}E_a\left(z,s\right)$$
for every $z\in H$, but we will use (6.3) during the proof
of Theorem 1.1, so it is enough for our purposes.

{\it Proof.\/} For $z\in H$ we have that
$$f\left(z\right)=\sum_{j=0}^{\infty}\left(f,u_j\right){\rm u}_j\left
(z\right)+\sum_a{1\over {4\pi}}\int_{-\infty}^{\infty}\left(f,E_a\left
(\ast ,{1\over 2}+ir\right)\right){\rm E}_a\left(z,{1\over 2}+ir\right
)dr,\eqno (6.4)$$
and this expression is uniformly and absolutely convergent
on compact subsets of $H$. Formula (6.4) follows from the Spectral Theorem
([I], Theorems 4.7 and 7.3).

Using Lemma 6.3 and that $f\left(\sigma_az\right)$ is also bounded, (6.4)
implies that for any $a\in A$ the sum
$$\sum_{j=0}^{\infty}\left(f,u_j\right)\beta_{a,j}\left(0\right)y^{
1-s_j}+$$
$$+\sum_c{1\over {4\pi}}\int_{-\infty}^{\infty}\left(f,E_c\left(\ast
,{1\over 2}+ir\right)\right)\left(\delta_{ca}y^{{1\over 2}+ir}+\phi_{
c,a}\left({1\over 2}+ir\right)y^{{1\over 2}-ir}\right)dr$$
is bounded for $z\in P\left(Y_{\Gamma}\right)$. Since by [I], (6.22) and (6.27) we have
$$\sum_c\left(f,E_c\left(\ast ,{1\over 2}+ir\right)\right)\phi_{c
,a}\left({1\over 2}+ir\right)=\left(f,E_a\left(\ast ,{1\over 2}-i
r\right)\right)$$
for any real $r$, so for any $a\in A$ the sum
$$\sum_{j=0}^{\infty}\left(f,u_j\right)\beta_{a,j}\left(0\right)y^{
1-s_j}+$$
$$+{1\over {4\pi}}\int_{-\infty}^{\infty}\left(\left(f,E_a\left(\ast
,{1\over 2}+ir\right)\right)y^{{1\over 2}+ir}+\left(f,E_a\left(\ast
,{1\over 2}-ir\right)\right)y^{{1\over 2}-ir}\right)dr$$
is bounded for $z\in P\left(Y_{\Gamma}\right)$, i.e.
$$\sum_{j=0}^{\infty}\left(f,u_j\right)\beta_{a,j}\left(0\right)y^{
1-s_j}+{1\over {2\pi i}}\int_{\left({1\over 2}\right)}\left(\int_
Ff\left(z\right)E_a\left(z,s\right)d\mu_z\right)y^{1-s}ds$$
is bounded for $z\in P\left(Y_{\Gamma}\right)$. We now shift
the integration to the right, to $\hbox{\rm Re$s=1-\delta$}$ with a small
$\delta >0$, and we see that every residue must be $0$ because of
the boundedness. The lemma follows.

 \bigskip\noindent {\bf References}

\nobreak
\parindent=12pt
\nobreak

\item{[B1]} A. Bir\'o, {\it On a generalization of the Selberg trace formula}, Acta Arithmetica, 87 (4) (1999), 319-338.

\item{[B2]} A. Bir\'o, {\it A relation between triple products of weight} $
0$ {\it and weight}  ${1\over 2}$ {\it cusp forms,\/} Israel J. of
Math., 182 (2011), no. 1, 61-101.

\item{[C]} F. Chamizo, {\it Some applications of large sieve in Riemann surfaces}, Acta Arithmetica, 77 (4) (1996), 315-337.

\item{[F]} J.D. Fay, {\it Fourier coefficients of the resolvent }
{\it for a Fuchsian group}, J. Reine Angew. Math., 294 (1977),
143-203.

\item{[G-R]} I.S. Gradshteyn, I.M. Ryzhik, {Table of integrals, series and
products, 6th edition,} {\it Academic Press}, 2000

\item{[He]} D.A. Hejhal, {The Selberg Trace Formula for $PSL(2,{\bf R})$, vol. 2, Lecture Notes in Mathematics 1001} {\it Springer, New York}, 1983

\item{[Hu]} H. Huber, {\it Ein Gitterpunktproblem in der hyperbolischen Ebene,\/} J. Reine Angew. Math., 496 (1998), 15-53.

\item{[H-J]} R.A. Horn, C.R. Johnson, {Matrix Analysis,} {\it Cambridge University Press Press}, 1985

\item{[I]} H. Iwaniec, {Introduction to the spectral theory of automorphic forms,} {\it Rev. Mat. Iberoamericana},
1995

\item{[K-N]} L. Kuipers, H. Niederreiter {Uniform distribution of sequences,} {\it New York, Wiley},
1974

\item{[R-P]} M. Risager, Y. Petridis, {\it Local average in hyperbolic lattice point counting, with an appendix by Niko Laaksonen}, Math. Z. (2016), doi:10.1007/s00209-016-1749-z, available online at https://arxiv.org/abs/1408.5743

\item{[Z1]} S. Zelditch, {\it Selberg trace formulae, pseudodifferential operators, and geodesic periods
of automorphic forms,} Duke Math. J., 56 (1988), no. 2, 295-344.

\item{[Z2]} S. Zelditch, {\it Trace formula for compact $\Gamma\setminus
PSL(2,{\bf R})$ and the equidistribution theory
of closed geodesics,} Duke Math. J., 59 (1989), no. 1, 27-81.

\item{[Z3]} S. Zelditch, {\it Selberg trace formulae and equidistribution theorems for closed geodesics
and Laplace eigenfunctions: finite area surfaces.,} Mem. Amer. Math. Soc. 96 (1992), no.
465, vi+102 pp

\bye